\documentclass[12pt]{elsart}
\usepackage{amsfonts,epsfig,latexsym}
\begin{document}


\thispagestyle{empty}

\begin{frontmatter}
\title{Small cones of $m$-hemimetrics}
\author{M.-M. Deza}
\address{CNRS/ENS, Paris, and Institute of Statistical Mathematics, Tokyo, {\tt Michel.Deza@ens.fr}}
\author{I. G. Rosenberg}
\address{D\'epartement de math\'ematiques et de statistique,
Universit\'e de Montr\'eal, {\tt rosenb@DMS.UMontreal.CA}}

\begin{abstract}

We introduce  polyhedral cones associated with 
$m$-hemimetrics on $n$ points, and, in particular, with $m$-hemimetrics coming
from partitions of an $n$-set into $m+1$ blocks. We compute generators and facets
of the cones for small values of $m,n$ and study their skeleton graphs.
\end{abstract}

\end{frontmatter}
\section{Introduction}

The notions of $m$-hemimetrics and $m$-partition hemimetrics are 
generalizations of the notions of metrics and cuts, which are 
well-known and central objects in Graph Theory, Combinatorial Optimization 
and, more generally, Discrete Mathematics.

Obviously the $m$-hemimetrics on an $n$-element set $E$
form a cone. What are its facets and extreme rays? The calculations
are rather complex due to the high dimension, $n \choose {m+1}$, of the cone
even for small $n$; moreover, matrices can no longer be used.

 Recall that a {\it metric} or a {\it metric space} is a pair $(E, d)$
where $E$ is a nonvoid set and $d: E^2\longrightarrow\mathbb{R}_+$(the
set of nonnegative reals) satisfies for all $x, y, z\in E$:
\\ (d1) $d(x, y) = 0 \iff x = y$,
\\(d2) $d(x,y)=d(y, x)$\quad\quad (symmetry),
\\(d3) $d(x,y)\leq d(x,z)+d(z,y)$\quad\quad (the triangle inequality).

A basic example is $(\mathbb{R}^2, d)$, where $d$ is the Euclidean
distance of $x$ and $y$; i.e., the length of the segment joining $x$
and $y$. An immediate extension is $(\mathbb{R}^3,d)$, where $d(x,y,z)$
is the area of the triangle with vertices $x, y$ and $z$. This leads
to the following definition (see  \cite{M,B,Fr,G1}). A
2-{\it metric} is a pair $(E, d)$, where $E$ is a nonempty set and $d:
E^3\longrightarrow\mathbb{R}_+$ satisfies for all $x, y, z, t\in E$
\\ (d1$^\prime$) $d(x, x,y) = 0$,
\\ (d1$^{\prime\prime}$) $x\neq y\Longrightarrow d(x,y,u)>0$ for some $u\in E$,
\\ (d2) $d(x,y,z)$ is totally symmetric,
\\ (d3) $d(x,y,z)\leq d(t,y,z)+d(x,t,z)+d(x,y,t)$ \ (the tetrahedron inequality).

The axiom (d2) means that the value of $d(x, y, z)$ is independent of
the order of $x, y$ and $z$. The axiom (d3) captures that fact that in
$\mathbb{R}^3$ the area of a triangle face of a tetrahedron does not
exceed the sum of the areas of the remaining three faces. 
A 2-metric allows the introduction of several geometrical and
topological concepts -- e.g. the betweenness, convexity, line and
neighborhood -- which lead to interesting results.

For finite 2-metrics (the object of this study), for their polyhedral aspects
and for applications, the axiom (d1$^\prime$) and  (d1$^{\prime\prime}$) seem
to be too restrictive and so we drop them. The definition given below
is formulated for an arbitrary positive integer $m$. A map $d:
E^{m+1}\longrightarrow\mathbb{R}$ is {\it totally symmetric} if for
all $x_1, ..., x_{m+1}\in E$ and every permutation $\pi$ of $\{1,..., m+1\}$$$
d(x_{\pi(1)}, ..., x_{\pi(m+1)}) = d(x_1, ..., x_{m+1}).$$

\noindent
{\bf Definition.} Let $m>0$. An {\it $m$-hemimetric} is a pair $(E, d)$,
where $d: E^{m+1}\longrightarrow\mathbb{R}$ is totally symmetric and
satisfies the {\it simplex inequality}: for all $x_1, ..., x_{m+2}$$\in E$$$
d(x_1, ..., x_{m+1}) \leq\sum^{m+1}_{i=1}d(x_1, ..., x_{i-1}, x_{i+1},
 ..., x_{m+2}).   \eqno (1)$$

Call a $m$-hemimetric $d$ {\it nonnegative} if $d$ takes only nonnegative 
values. The notion of a $m$-hemimetric is new, but it is closely related to
the notion of a $m$-metric, considered in about 200 references collected in \cite{G2}.
The study of $m$-hemimetrics is motivated (apart, of course, that it
represents an extension of metrics) by applications in Statistics and Data
Analysis (see \cite{DR} for some relevant references). Notice the following immediate:

\noindent
{\bf Fact 1.} 
{\it If $(E, d)$ and $(E, d^\prime)$ are $m$-hemimetrics and $a,
b\in\mathbb{R}_+$ then $(E, ad+bd^\prime)$ is an $m$-hemimetric.}

Here, as usual, for all $x_1, ..., x_{m+1}\in E$$$
(ad+bd^\prime)(x_1, ...,x_{m+1}) := ad(x_1, ..., x_{m+1})+
bd^\prime(x_1, ..., x_{m+1}).$$

Given an $(m+1)$-partition  $S_{1}, \ldots, S_{m+1}$ of 
$V_n: = \{1, 2, \ldots, n\}$, a {\it partition $m$-hemimetric}
$\alpha(S_{1}, \ldots, S_{m+1})$ is defined by setting
$\alpha(S_{1}, \ldots, S_{m+1}) (i_1, \dots, i_{m+1})$ is equal to 1 if 
for no $1 \le j<l \le m+1$ both $i_j$ and $i_l$ belong to the same $S_k$, and 0 otherwise.
It is easy to see that $\alpha(S_{1}, \ldots, S_{m+1})$ is a nonnegative $m$-hemimetric
and for $m=1$ it is the usual {\it cut} semimetric (see, for example,
\cite{DL}). 

For small values of $n$ and $m$ we consider the cone of all $m$-hemimetrics,
the cone of all nonnegative  $m$-hemimetrics and the cone,
generated by all partition $m$-hemimetrics on $V_n$.
 Using computer search we list facets and generators for 
these cones and tables of their adjacencies and incidences. The different
orbits were determined manually, using symmetries.
We study two graphs, the $1$-skeleton and the ridge graph, of these 
polyhedra: the number of their nodes and edges,
 their diameters, conditions of adjacency, inclusions among the graphs 
and their restrictions on some orbits of nodes.
In fact, we would like to describe two graphs $G(C),  G(C^*)$ for our three cones
as fully as possible, but in the cases, when it is too difficult, we will give some partial
information on adjacencies in those graphs. Especially we are interested in
the diameters of the graphs, in a good criterion of adjacency, in their 
local graphs
(i.e. in the subgraphs induced by all neighbors of a given vertex) and in their
restrictions on some orbits. Finally, we compare obtained
 results with similar results for metric case (see
 ~\cite{DD2,DDF,DL}) and quasi-metric case (see
\cite{DP}).
 All  computation was done using the programs {\it cdd} of ~\cite{Fu}.
\newline

The following  notation will be used below:
\begin{itemize}
\item the {\it (m+1)-simplex} inequality (1) and, in particular, for $m=2$,
the {\it tetrahedron} inequality

$T_{ijk, l}: x_{ijl} + x_{ikl} + x_{jkl} - x_{ijk} \geq 0$;
\item the {\it nonnegativity} inequality
$N_{i_1,\dots,i_{m+1}}: x_{i_1,\dots,i_{m+1}} \geq 0$;  
\item the {\it cone $P_{n}^m$ of partition m-hemimetrics}, generated by all 
$(m+1)$-partitions of $V_n$;
\item the {\it  cone $NHM_{n}^m$ of nonnegative m-hemimetrics}, 
defined by all

$(n-m-1)$$n \choose {m+1}$ \hspace{3mm} $(m+1)$-simplex
inequalities and all $n \choose {m+1}$ nonnegativity inequalities on $V_n$;    
\item the {\it  cone of all m-hemimetrics}, $HM_{n}^m$, defined by 
all $(m+1)$-simplex inequalities on $V_n$.  
\end{itemize}

Clearly, $P_n^m \subseteq NHM_n^m \subseteq HM_n^m$ and these 3 cones are of
full dimension $n \choose {m+1}$ each. For $m=1$ the last two cones coincide
and
the first two cones are the {\it cut cone $CUT_n$} and the {\it semimetric cone
$MET_n$}, considered in detail in \cite{DL} and in the references listed
there.
The cone $P_n^1=CUT_n$ has $ \lfloor \frac{n}{2} \rfloor$ orbits of extreme
rays
(represented by the cuts $ \alpha(1 \dots i,(i+1) \dots n)$).

To simplify the notation we keep $n$ fixed and denote by $E_k$ the family of all 
$k$-element subsets of $V_{n}$ ($k=1,...,n$). Let $d$ be a semimetric on the set $V_{n}$. 
Because of symmetry (d2) and since $d(i,i) =0$ for all $i \in V_{n}$, we can view 
the semimetric $d$ as a vector $(d_{i_1,i_2})_ \in 
\mathbb{R}^{E_{2}}$. In the same way, we 
can view the $m$-hemimetric $d$ on the set $V_{n}$  as a vector 
$(d_{i_1,\dots,i_{m+1}})_ \in \mathbb{R}^{E_{m+1}}$. 
In particular, each extreme ray of the cones, considered below, will be
represented by an integer vector on the ray with relatively prime coordinates.
We also will
represent the facets of cones by such vectors.

Any such vector $v=(v_{i_1,\dots,i_{m+1}}) \in \mathbb{R}^{E_{m+1}}$ can be represented
by the following vertex-labeled induced subgraph $R(v)$ of the Johnson graph
$J(n,m+1)$. The vertices of $R(v)$ are all
unordered $(m+1)$-tuples $(i_1,\dots,i_{m+1})$ such that
$(v_{i_1,\dots,i_{m+1}})$ is not
zero. This value will be the label of the vertex; we will omit the label when
it is 1. Two vertices of the Johnson graph (and also of its induced subgraph 
$R(v)$) are adjacent if the corresponding $(m+1)$-tuples have $m$ common elements. 

For example,
$R( \alpha(S_1 , \dots, S_{m+1}))$ is the complement to the Hamming graph
$H( \vert S_{1} \vert, \ldots,\vert S_{m+1} \vert)$, i.e. the direct (Cartesian) product of the cliques $K_{ \vert S_i \vert }$, $1 \le i \le {m+1}$.
Another example: the graph $R$ of the vector defining a nonnegativity facet is a vertex and, for a
$(m+1)$-simplex facet, it is the complete graph $K_{m+2}$ with one vertex labeled $-1$.

\section{Partition $m$-hemimetrics and related polyhedra}

Recall that for a partition $S_1,S_2$ of $V_n$ the cut semimetric
$\alpha(S_1, S_2)$ satisfies  $\alpha(S_1, S_2)_{ij}=1$ if
 $\{i,j\} \cap S_1$ is a singleton and
$\alpha(S_1, S_2)_{ij}=0$ otherwise. We extend it as follows.
Let $q \geq 2$ be an 
integer and let $S_{1}, \ldots, S_{q}$ be pairwise disjoint nonvoid subsets of 
$V_{n}$, forming a partition of $V_{n}$. 
The {\it multicut semimetric} $\delta (S_{1}, \ldots, 
S_{q})$ is the vector in  $\mathbb{R}^{E_{2}}$, defined by $\delta (S_{1}, \ldots, 
S_{q})_{ij} = 0$, if $i, j \in S_{h}$ for some $h$, $1 \leq h \leq q$, and 
$\delta (S_{1}, \ldots, S_{q})_{ij} = 1$, otherwise.

The connection between $\delta (S_{1}, \ldots,S_{q})$ and
$\alpha (S_{1}, \ldots, S_{q})$ from Section 1 is given by

 $\alpha (S_{1}, \ldots, S_{q})(i_1,\dots,i_{q})$=$ \prod_{1 \le s<t \le q}\delta (S_{1}, \ldots,S_{q})(i_s,i_t)$=

$ \lfloor \frac {\sum_{1 \le s<t \le q} \delta (S_{1},\ldots,S_{q})(i_s,i_t)}{{q \choose 2}} \rfloor$; compare it with the half-perimeter $m$-semimetric from \cite{DR}.

The cone generated by all multicut semimetrics
$\delta(S_{1}, \ldots, S_{q})$ ($q \ge 2$) on $V_{n}$, is called {\it the multicut cone}
and  denoted by $MCUT_{n}$; it coincides with $CUT_{n}$ (see \cite{DL},
Proposition 4.2.9).
The convex hull of the cut semimetrics (multicut semimetrics) on 
$V_{n}$, is called the {\it cut polytope} ( {\it multicut polytope}) and is
denoted by $CUT^{\Box}_{n}$ ($MCUT^{\Box}_{n}$); the two polytopes not coincide.

\section{ Facets, extreme rays and their orbits in polyhedra}

We recall some terminology. Let $C$ be a polyhedral cone in $\mathbb{R}^{n}$. Given $v \in \mathbb{R}^{n}$, the 
inequality $v^{T}x \leq 0$ is said to be {\it valid} for $C$, if it 
 holds for all $x \in C$. Then the set $\{ x \in C \vert v^{T}x = 0 \}$ 
is called the {\it face of $C$, induced by the valid inequality $v^{T}x \leq 
0$}. A face of dimension $\dim(C) - 1$ is called a {\it facet} of $C$;
 a face of dimension $1$ is called an {\it extreme ray} of $C$.
A face of dimension  $\dim(C) - 2$ is called a {\it ridge}.

Two vertices $x, y$ of $C$ are said to be {\it adjacent}, if they
generate a face of dimension 2 of
$C$. Two facets of $C$ are said to be {\it adjacent}, 
if their intersection has dimension $\dim(C) - 2$. 
The {\it $1$-skeleton} graph of $C$ is the graph 
$G(C)$ whose nodes are the extreme rays of $C$ and whose edges
are the pairs of adjacent nodes. Denote by $C^*$ the dual cone of $C$. The
{\it ridge graph} of $C$ is the graph whose nodes are
the facets of $C$ and with an edge between two facets if
they are adjacent on $C$. So, the ridge graph of a cone $C$ is
the $1$-skeleton $G(C^*)$ of its dual cone.
 
A mapping $f: \mathbb{R}^{n} \longrightarrow \mathbb{R}^{n}$ is called a {\it symmetry} of  
a cone $C$ (or a polytope $P$), if it is an isometry, satisfying $f(C)=C$
(or $f(P)=P$). ({\it An isometry} of $\mathbb{R}^{n}$ is a linear mapping preserving
the Euclidean distance.)
Given a face $F$,  the {\it orbit} $\Omega(F)$ of $F$ consists of 
all  faces, that can be obtained from $F$ by the group of all symmetries of $C$.
 
Clearly, all the faces of $CUT_{n}$ and $CUT^{\Box}_{n}$ are preserved by any
permutation of $V_n$.

For $m>1$ all orbits of faces of $m$-hemimetric cones
$P_n^m, NHM_n^m, HM_n^m$ on $V_{n}$ 
are also preserved under any permutation of the set $V_{n} = \{1, \ldots, n\}$.
We conjecture that the symmetry group  
consists only of permutations of $V_n$, i.e. it is the group $Sym(n)$ of all
permutations on $V_n$ ( see Theorem 3.3 in \cite{DGL} stating that the symmetry group of
a truncated multicut polytope is $Sym(n)$).

 \section{The case of $n=m+2$}

The minimal $n$ for which the three cones are nontrivial is $m+2$; the dimension of the cones is also $m+2$ for
$n=m+2$.

 First, we present a complete linear description  for case
$(m,n) =(2,4)$.

It turns out that $P_4^2=NHM_4^2$. This cone has 6 extreme rays (all in
the same orbit under $Sym(4)$): $ \alpha(S_1,S_2,S_3)$ for the 3-partitions

$(1,2,34), (1,3,24), (1,4,23), (2,3,14), (2,4,13), (3,4,12)$.

There are $8$ facets, which form $2$ orbits: the orbit $F_{1}$ of all 4
tetrahedron facets and the orbit $F_{2}$ of all 4 nonnegativity facets.

The edge graph $G(C)$ is $K_6-3K_2$ (the octahedron); the 3 pairs of nonadjacent
rays are of the form $\alpha(a,b,cd), \alpha(c,d,ab)$. Each extreme ray (say, $\alpha(1,2,34)$)
is incident to 2 tetrahedron and to 2 nonnegativity facets (namely, to
$T_{123,4}, T_{124,3}$ and $N_{134}, N_{234}$).

The ridge graph $G(C^*)$ is the cube.
Adjacencies of facets of $NHM_4^2$ are shown in Table 
\ref{tab:tabl1}. For each orbit a representative and the number of adjacent 
facets from other orbits are given, as well as the total
 number of adjacent ones, the number of incident extreme rays and the cardinality of orbits. 

\begin{table}
\caption{The adjacencies of facets in the cone $NHM_4^2$}
 \begin{tabular}{|c|c|c|c|c|c|c|} 
 \hline 
 Orbit & Representative & $F_{1}$ &  $F_{2}$ & Adj.&Inc. & $\vert F_{i} 
 \vert $ \\ 
 \hline 
 $F_{1}$ & $T_{123, 4}$ & $0$ & $3$ & $3$&$3$ & $4$ \\ 
 $F_{2}$ & $N_{123}$ & $3$ & $0$ & $3$&$3$  & $4$ \\ 
 \hline 
 \end{tabular}
\label{tab:tabl1}
\end{table}

More precisely, for the ridge graph of $NHM_4^2$ it holds:
    
   (i) The tetrahedron facet $T_{ijk,l}$ is adjacent only to the
facets $N_{ijl}, N_{ikl}, N_{jkl}$;

   (ii) The nonnegativity facet $N_{ijk}$ is adjacent only to the facets 
$T_{ijl,k}, T_{ikl,j}, T_{jkl,i}$.

The cone $HM_{m+2}^m$ is a simplex $(m+2)$-dimensional cone; so
$G(C)=G(C^*)=K_{m+2}$. Its facets are all $(1,-1)$-valued
$(m+2)$-vectors with only one $-1$, its generators are all $(1-m,1)$-valued $(m+2)$-vectors
with only one $1-m$. Notice that
$P_3^1=HM_3^1=CUT_3=MET_3$.

In general, $P_{m+2}^m=NHM_{m+2}^m$ for any $m \ge 2$. This cone 
has ${m+2} \choose 2$ extreme rays, all in the same orbit, represented
by $ \alpha(12,3, \dots, m+2)$, i.e. by any vector of length $m+2$, consisting
of two ones and $m$ zeros. The skeleton of $P_{m+2}^m$ is the Johnson graph
$J(m+2,2)$, called also the
{\it triangular} graph $T(m+2)$, which is is the line graph $L(K_{m+2})$. It is also
the skeleton of the $(m+1)$-polytope (called ${\it ambo}- \alpha_{m+1}$), obtained from the
$(m+1)$-simplex as the convex hull of the mid-points of all its edges; e.g.
$T(4)$ is the skeleton of the octahedron, $T(5)$ is the complement of the Petersen
graph. In general, $T(m), m \ge 2,$ has diameter 2; moreover, it is a strongly regular graph.

The cone $P_{m+2}^m$ has two orbits, $F_1$ and $F_2$, of facets, 
containing $m+2$ facets each and represented by the $(m+1)$-simplex facet
$T_{1\dots(m+1),(m+2)}$ and by the nonnegativity facet $N_{1\dots(m+1)}$. The
orbit $F_1$ consists of simplex cones, i.e. facets from this orbit are incident to $m+1$ linearly
independent extreme rays. Any nonnegativity inequality $N$  defines
the cone $P_{m+1}^{m-1}=NHM_{m+1}^{m-1}$, i.e. it becomes equality on this
smaller cone. So, $N$ is non-facet only for $m=1$ and it is a simplex cone
only for $m=2$; in general, $N$ is incident to ${m+1} \choose 2$ extreme rays.
The ridge graph is $ \overline {K_{m+2}}$ on $F_1$; on $F_2$ it is
$ \overline {K_4}$ for $m=2$ and $K_{m+2}$ for $m \ge 3$. Finally, the $m+2$
pairs $(T_{ \overline i, i}, N_{ \overline i})$ (of $(m+1)$-simplex and
nonnegativity facets) are the only non-edges for pairs of
facets from different orbits.

\section{Small $2$-hemimetrics}

\subsection {The case of $5$ points}

We present here the complete linear description of $P_5^2$, $NHM_5^2$ and
$HM_5^2$. The cone  $P_5^2$ has $25$ extreme rays, which form $2$ orbits
with representatives $\alpha(1,2,345)$ (orbit $O_1$) and
$\alpha(1,23,45)$ (orbit $O_2$). The skeleton and the ridge graph of $P_5^2$
has 270 and 1185 edges, respectively. The cone $P_5^2$ has $120$ facets divided
into $4$ orbits, induced by the $20$ tetrahedron inequalities (orbit $F_{1}$),
the $10$ nonnegativity inequalities (orbit $F_{2}$),
the $60$ inequalities (orbit $F_{3}$), represented by 

 \[A: 2x_{123} - (x_{124} + x_{135}) + (x_{134} + x_{125} + x_{245} 
 +x_{345}) \ge 0\] 
 and the $30$ inequalities (orbit $F_{4}$), represented by

 \[B: 2(x_{123} + x_{145} + x_{245} - x_{345}) + (x_{134} + x_{135} +
 x_{234} + x_{235} - x_{124} - x_{125}) \ge 0\].

The above two inequalities are the 2-hemimetric analogs of the
following {\it 5-gonal}
inequality (the simplest inequality, different from
the triangle inequality), appearing  in the cone $CUT_n$ for $n \ge 5$:

\[ (x_{13} + x_{14} + x_{15} + x_{23} + x_{24} + x_{25}) -
 (x_{12} + x_{34} + x_{35} + x_{45}) \ge 0\].

This facet and $B$ have both the Petersen graph as their $ \overline R$ (i.e. 
the complement of their graph $R$).
Clearly, the graphs
$R$ for partition 2-hemimetrics $\alpha(1,2,34)$, $\alpha(1,2,345)$,
$\alpha(1,23,45)$ are the cycles $C_2, C_3, C_4$. 

The graphs $ \overline {R(A)}$, $ \overline {R(B)}$ are given on the Figure 1.

\begin{center}
\epsfysize=4cm
\epsffile{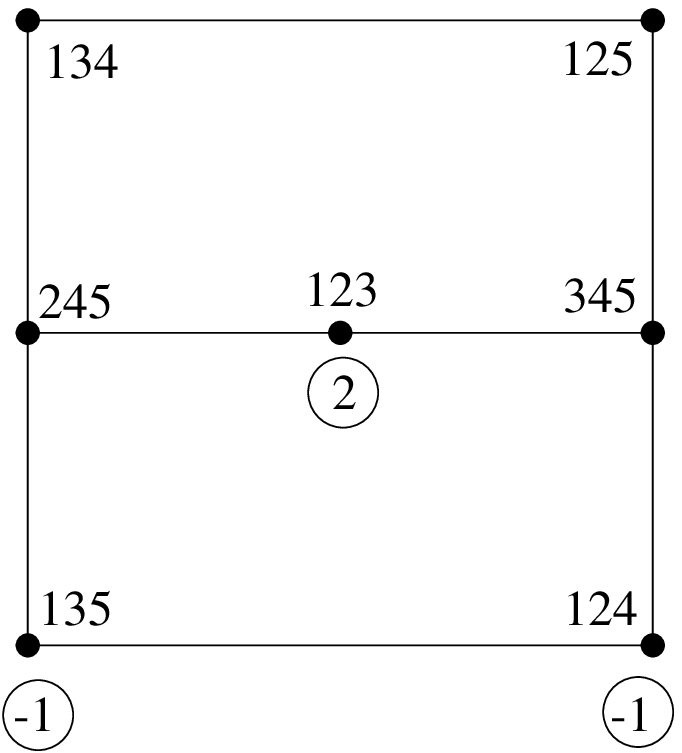}
\hspace{1cm}
\epsfysize=4cm
\epsffile{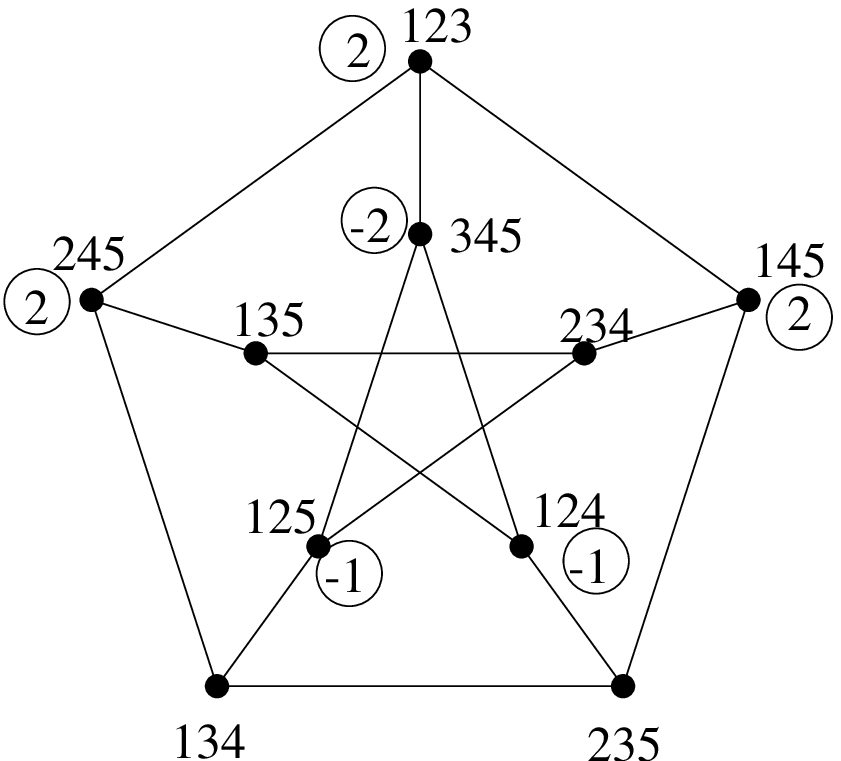}
\end{center}
\begin{center}
Figure 1: $ \overline {R(A)}$, $ \overline {R(B)}$ in the cone $P_5^2 $
\end{center}

The facets from the orbit $F_4$ are simplex-cones, i.e. the extreme rays
on them are linearly independent. Among the 9 neighbors of $B$, 4 are from the
orbit $F_1$, 4 are from the orbit $F_3$ and exactly one (actually,
$N_{123}$) from the orbit $F_2$. The local graph of a facet from $F_4$
(i.e. the subgraph of the ridge graph of $P_5^2$, induced by all neighbors
of $B$) is $K_9 - C_4$. In fact, all nonadjacencies in this local graph are
the four edges of the 4-cycle of the 4 neighbors of $B$ from the orbit $F_3$.

Facets from orbits $F_1,F_2,F_3,F_4$ are incident, respectively, to 7,9; 7,9;
 4,6; 3,6 extreme rays from orbits $O_1,O_2$ of $P_5^2$.

 \begin{table}
\caption{The adjacencies of extreme rays in the cone $P_5^2$}
 \begin{tabular}{|c|c|c|c|c|c|c|} 
 \hline 
 Orbit & Representative & $O_{1}$ &  $O_{2}$ & Adj.&Inc. & $\vert O_{i} 
 \vert $ \\ 
 \hline 
 $O_{1}$ & $\alpha(1,2,345)$ & $9$ & $12$ & $21$&$54$ & $10$ \\ 
 $O_{2}$ & $\alpha(1,23,45)$ & $8$ & $14$ & $22$&$54$  & $15$ \\ 
 \hline 
 \end{tabular}
\label{tab:tabl2}
\end{table}

 \begin{table}
\caption{The adjacencies of facets in the cone $P_5^2$}
 \begin{tabular}{|c|c|c|c|c|c|c|c|c|} 
 \hline 
 Orbit & Representative & $F_{1}$ &  $F_{2}$ &$F_{3}$&$F_{4}$ & Adj.&Inc. & $\vert F_{i} 
 \vert $ \\ 
 \hline 
 $F_{1}$ & $T_{123, 4}$ & $16$ & $9$ & $18$&$6$ & $49$&$16$&$20$ \\ 
 $F_{2}$ & $N_{123}$ & $18$ & $3$ & $18$&$3$  & $42$&$16$&$10$ \\
 $F_{3}$ & $A$ & $6$ & $3$ & $4$&$2$  & $15$&$10$&$60$ \\
 $F_{4}$ & $B$ & $4$ & $1$ & $4$&$0$  & $9$&$9$&$30$ \\
 \hline 
 \end{tabular}
\label{tab:tabl3}
\end{table}

The skeleton and the ridge graph of $NHM_5^2$ have 420 and 355 edges, respectively. The adjacencies of 37 extreme rays and of 30 facets of this cone are
given in Tables \ref{tab:tabl4}, \ref{tab:tabl5}. The extreme rays are divided
into 3 orbits $O_1,O_2,O_3$, represented by (0,1)-valued vectors $v_1,v_2,v_3$
below; their $R$-graphs are $C_3,C_4,C_5$, respectively.

 \begin{table}
\caption{The adjacencies of extreme rays in the cone $NHM_5^2$}
 \begin{tabular}{|c|c|c|c|c|c|c|c|} 
 \hline 
 Orbit & Representative & $O_{1}$ &$O_{2}$&$O_{3}$& Adj.&Inc. & $\vert O_{i} 
 \vert $ \\ 
 \hline 
 $O_{1}$ & $\alpha(1,2,345)$ & $9$ & $12$ & $6$&$27$&$21$ & $10$ \\ 
 $O_{2}$ & $\alpha(1,23,45)$ & $8$ & $6$ & $8$&$22$&$18$  & $15$ \\
 $O_{3}$& $v_3$&$5$&$10$&$5$&$20$&$15$& $12$ \\
 \hline 
 \end{tabular}
\label{tab:tabl4}
\end{table}

 \begin{table}
\caption{The adjacencies of facets in the cone $NHM_5^2$}
 \begin{tabular}{|c|c|c|c|c|c|c|} 
 \hline 
 Orbit & Representative & $F_{1}$ &  $F_{2}$ & Adj.&Inc. & $\vert F_{i} 
 \vert $ \\ 
 \hline 
 $F_{1}$ & $T_{123, 4}$ & $16$ & $9$ & $25$&$22$ & $20$ \\ 
 $F_{2}$ & $N_{123}$ & $18$ & $3$ & $21$&$22$  & $10$ \\ 
 \hline 
 \end{tabular}
\label{tab:tabl5}
\end{table}

Each facet (from both orbits) of $NHM_5^2$ is incident to 7,9,6 extreme rays
from orbits $O_1,O_2,O_3$, respectively. Each (tetrahedron) facet of $HM_5^2$
 is incident to 7,9,6,6,9,15 extreme rays from orbits $O_1,\dots,O_6$.

For any cone, let $I_{O_i,F_j}$ and $ I_{F_j,O_i}$ denote the number of facets from the orbit $F_j$,
incident to an extreme ray of the orbit $O_i$, and, respectively, the number
of extreme rays from $O_i$, incident to a facet from $F_j$. Clearly,
$\vert O_{i}\vert I_{O_i,F_j} = \vert F_{j}\vert I_{F_j,O_i} $.       

The cone $HM_5^2$ has 92 extreme rays divided into 6 orbits. Below we give
some representatives $v_1,\dots,v_6$ of those orbits $O_1,\dots,O_6$.
The first two represent both orbits of $P_5^2$, the first 3 represent the
3 orbits of $NHM_5^2$.

  $x=(x_{123}, x_{124}, x_{125}, x_{134}, x_{135}, x_{145}, x_{234}, x_{235}, x_{245}, x_{345})$:

$v_1=(1,     1,     1,      0,     0,     0,     0,     0,     0,     0    )$;

$v_2=(0,     1,     1,      1,     1,     0,     0,     0,     0,     0    )$;

$v_3=(1,     0,     1,      0,     0,     1,     1,     0,     0,     1    )$;

$v_4=(1,     1,     1,     -1,     0,     0,     1,     0,     1,     1    )$;

$v_5=(1,     1,     1,     -1,    -1,     1,     1,     1,     1,     1    )$;
 
$v_6=(1,     0,     1,      0,     1,    -1,     1,     1,     2,     1    )$.

 \begin{table}
\caption{The adjacencies of extreme rays in the cone $HM_5^2$}
 \begin{tabular}{|c|c|c|c|c|c|c|c|c|c|c|} \hline Orbit & Representative & 
 $O_{1}$ & $O_{2}$ & $O_{3}$ & $O_{4}$ & $O_{5}$&$O_6$ &Adj.&Inc. & $\vert 
 O_{i} \vert $ \\ \hline
 $O_{1}$ &$v_1=\alpha(1,2,345)$&$6$&$12$&$6$&$6$& $9$ & $9$ & $48$ & $14$ & $10$ \\
 $O_{2}$ &$v_2=\alpha(1,23,45)$&$8$&$2$&  $4$&$2$&$4$&$8$&$28$&$12$&$15$\\
 $O_{3}$ &$v_3$ &$5$&$5$& $0$ & $5$ & $5$ & $5$ & $25$ & $10$ & $12$ \\
 $O_{4}$ &$v_4$ & $6$&$3$ & $6$&$0$ & $3$ & $6$ & $24$ & $12$ & $10$ \\
 $O_{5}$ &$v_5$ & $6$&$4$&$4$ & $2$ & $0$ & $4$ & $20$ & $12$ & $15$ \\
 $O_{6}$ &$v_6$ &$3$& $4$ & $2$ &$2$& $2$ & $0$ & $13$ & $10$ & $30$ \\
 \hline 
\end{tabular}
\label{tab:tabl16}
 \end{table}

\begin{prop}

 The diameters of the skeleton graphs of $P_5^2$ and of $NHM_5^2$ are $2$.
 \end{prop}

In fact, each of the orbits $O_1, O_2$ of $P_5^2$ is a dominating clique. There is only one 
type of a non-edge, represented by $\alpha(1,23,45), \alpha(2,3,145))$, but
$\alpha(1,3,245)$ is one of common neighbors. The complement of the skeleton
of $P_5^2$ turns out to be the Petersen graph with a new vertex (corresponding
to a member of the orbit $O_2$) on each of 15 edges. The result for $NHM_5^2$
comes also by finding out a common neighbor to each possible non-edge.

\begin{prop} For the ridge graphs of $NHM_5^2$ and $HM_5^2$ it holds:   

(i) The diameter of the ridge graph of $NHM_5^2$ is $2$;

(ii) Its restriction on the orbits $F_1$ and $F_2$ is $K_{4,4,4,4,4}$ and the
Petersen graph, respectively;

(iii) The ridge graph of $HM_5^2$ is $K_{4,4,4,4,4}$ (of diameter $2$).
\end{prop}

\subsection{The case of $6$ points}

$NHM_6^2$ has exactly 12492 extreme rays, with
$(adjacency, incidence)$ pairs being, respectively, (2278,64), (1321,56), (1030,40), (818,48), (731,48), (358,40), (270,36), (93,28), (66,28), (51,28), (47,28), (46,39), (37,31), (32,28), (30,27), (29,26), (27,23), (26,24), (26,23), (25,25), (23,22), (22,21), (21,21).

Three of the above pairs (1st, 2nd and 4th) are realized by orbits (say,
$0_1$, $O_2$ and $O_4$), which are represented by 3-partition
2-hemimetrics $\alpha(1,2,3456)$, $\alpha(1,23,456)$, $\alpha(12,34,56)$
and have size 15, 60, 15, respectively. The graphs $R$ of members
of orbits $O_1, O_2$ and  $O_4$ are $K_4, K_6-C_6 = K_3 \times K_2$ (the skeletons
of the tetrahedron and 3-prism) and the skeleton of the cube. 
Three other orbits consist also of
(0,1)-valued extreme rays: $O_3$ (with $R$ being  the Petersen graph),
$O_5$ (with $R$-graph being the skeleton of the simple polyhedron with $p$-vector
$p=(p_3=2,p_4=2,p_5=2)$) and $O_7$ with graph $R$ (non-planar,
 non-regular), given on Figure 2, together with one for $O_5$.
The extreme rays of the remaining orbits are (0,1,2)-valued and
(0,1,2,3)-valued vectors.

\begin{center}
\epsfysize=4cm
\epsffile{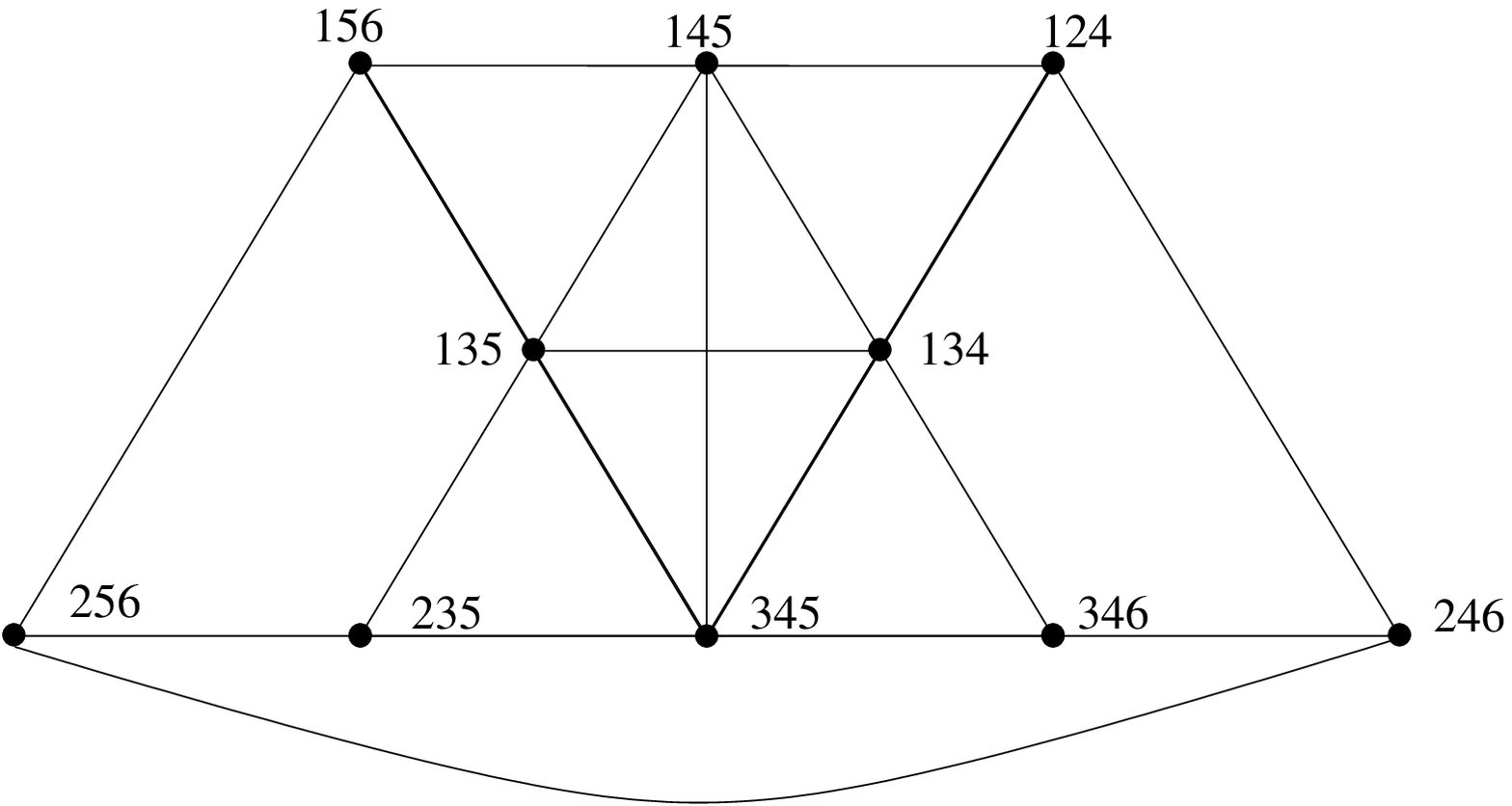}
\hspace{1cm}
\epsfysize=4cm
\epsffile{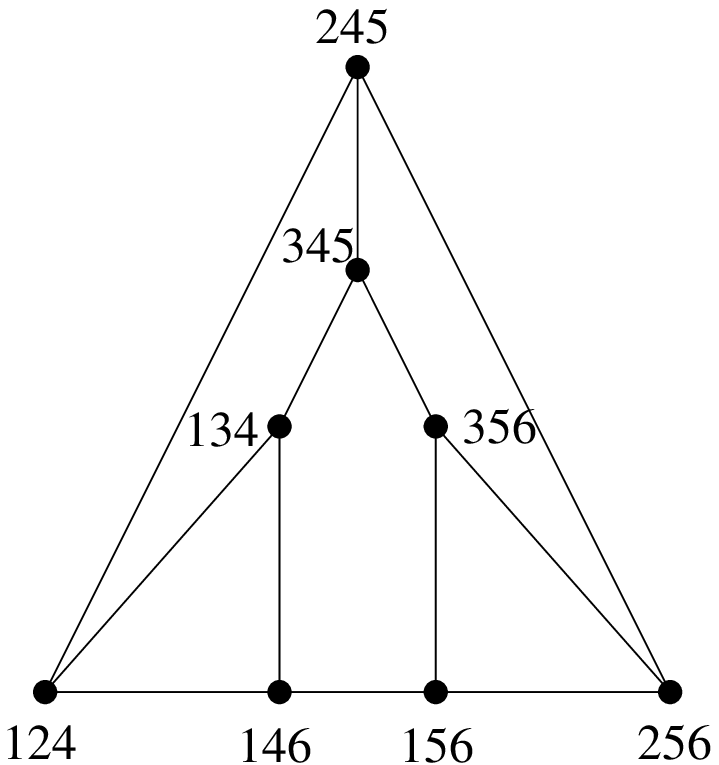}
\end{center}
\begin{center}
Figure 2 : The graphs R of extreme rays from orbits $ O_7, O_5 $ of the cone 
$NHM_6^2$
\end{center}

\begin{table}
\caption{The adjacencies of facets in the cone $NHM_6^2$}
 \begin{tabular}{|c|c|c|c|c|c|c|} 
 \hline 
 Orbit & Representative & $F_{1}$ &$F_{2}$& Adj.&Inc. & $\vert F_{i} \vert  $ \\
\hline
 $F_{1}$ & $T_{123,4}$ & $56$ & $19$ & $75$&$4001$  & $60$ \\ 
 $F_{2}$ & $N_{123}$ & $57$ & $10$ & $67$&$3939$  & $20$ \\ 
 \hline 
 \end{tabular}
\label{tab:tabl7}
\end{table}

The cone $P_6^2$ has more than 950.000 facets (computer stopped,
by lack of memory, after 72, out of 90, iterations).
Here are two examples of a $(0,1,-1)$-valued facets of $P_6^2$;
see also Figure 3 (for the facet W).

$W: (-x_{145}+x_{146}+x_{136}+x_{123}+x_{125}) + (x_{245}+x_{234}+x_{346}+x_{356}+x_{256}) - (x_{235}+x_{236}) \ge 0.$

\begin{center}
\epsfysize=4cm
\epsffile{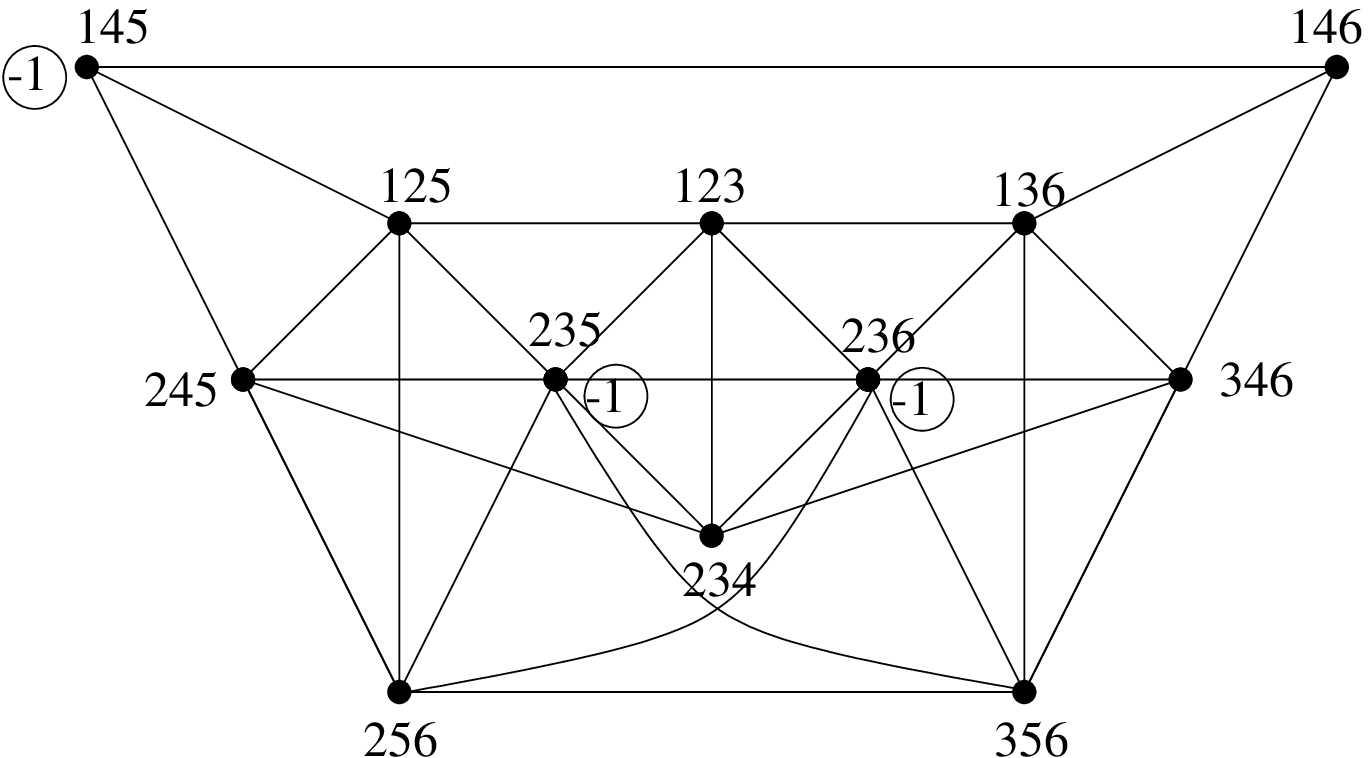}
\end{center}
\begin{center}
Figure 3 : $ R(W) $ in the cone $P_6^2$
\end{center}

$Z: \sum x_{ijk}-(x_{124}+x_{125}+x_{145})-(x_{234}+x_{235}+x_{345})-
2(x_{146}+x_{156}+x_{456})-2x_{236} \ge 0.$

Remark that the triples with coefficients zero, in $W$ and $Z$, form
the skeleton of 1- and 2-truncated tetrahedron, respectively; the triples
with coefficient -1 form $K_3+K_1$ and $K_2+K_1$, respectively.

\section{Small $3$-hemimetrics}

The cone $NHM_6^3$ has 287 extreme rays divided into 5 orbits. Below we give
representatives $u_1,\dots,u_5$ of the orbits $O_1,\dots,O_5$.
These vectors are indexed by 4-subsets of the set $\{1,\dots,6\}$; the 4-subsets are given as the complements of 2-subsets. The first four are $(0,1)$-valued; their $R$-graphs (in the Johnson
graph $J(6,4)$ of all $4$-tuples) are the cycles $C_3, C_4, C_5, C_6$, respectively.
The first two are partition $3$-hemimetrics; they represent both orbits
of $P_6^3$. The graphs $R(u_4)$ and $R(u_5)$ are on Figure 4.

  $x=(x_{\overline{12}}, x_{\overline{13}}, x_{\overline{14}}, x_{\overline{15}}, x_{\overline{16}}, x_{\overline{23}}, x_{\overline{24}}, x_{\overline{25}},
 x_{\overline{26}}, x_{\overline{34}}, x_{\overline{35}}, x_{\overline{36}},
x_{\overline{45}}, x_{\overline{46}}, x_{\overline{56}})$:

$u_1=(0,  0,  1,  1,  0,  0,  0,  0,  0,  0,  0,  0,  1,  0,  0)$;

$u_2=(0,  0,  1,  1,  0,  0,  0,  0,  0,  1,  1,  0,  0,  0,  0)$;

$u_3=(0,  0,  1,  1,  0,  0,  0,  0,  0,  1,  0,  1,  0,  0,  1)$;

$u_4=(0,  0,  1,  1,  0,  0,  0,  1,  1,  1,  0,  1,  0,  0,  0)$;

$u_5=(0,  0,  1,  1,  0,  1,  0,  0,  1,  0,  0,  1,  1,  2,  0)$;

\begin{center}
\epsfysize=4cm
\epsffile{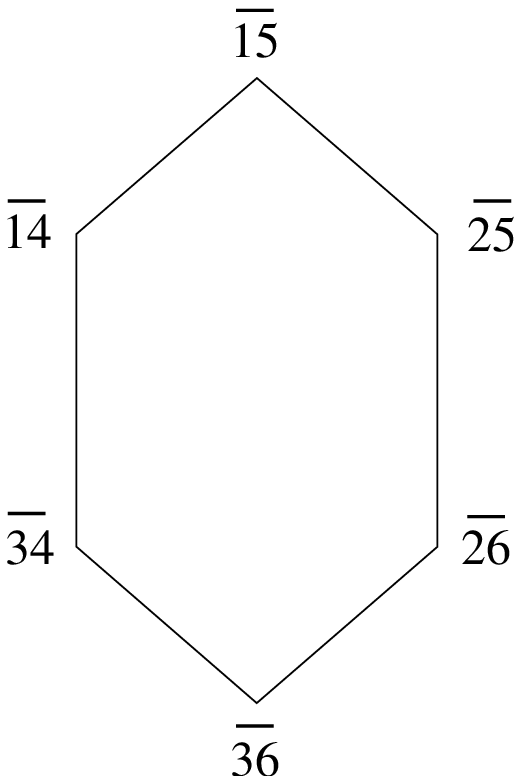}
\hspace{1cm}
\epsfysize=4cm
\epsffile{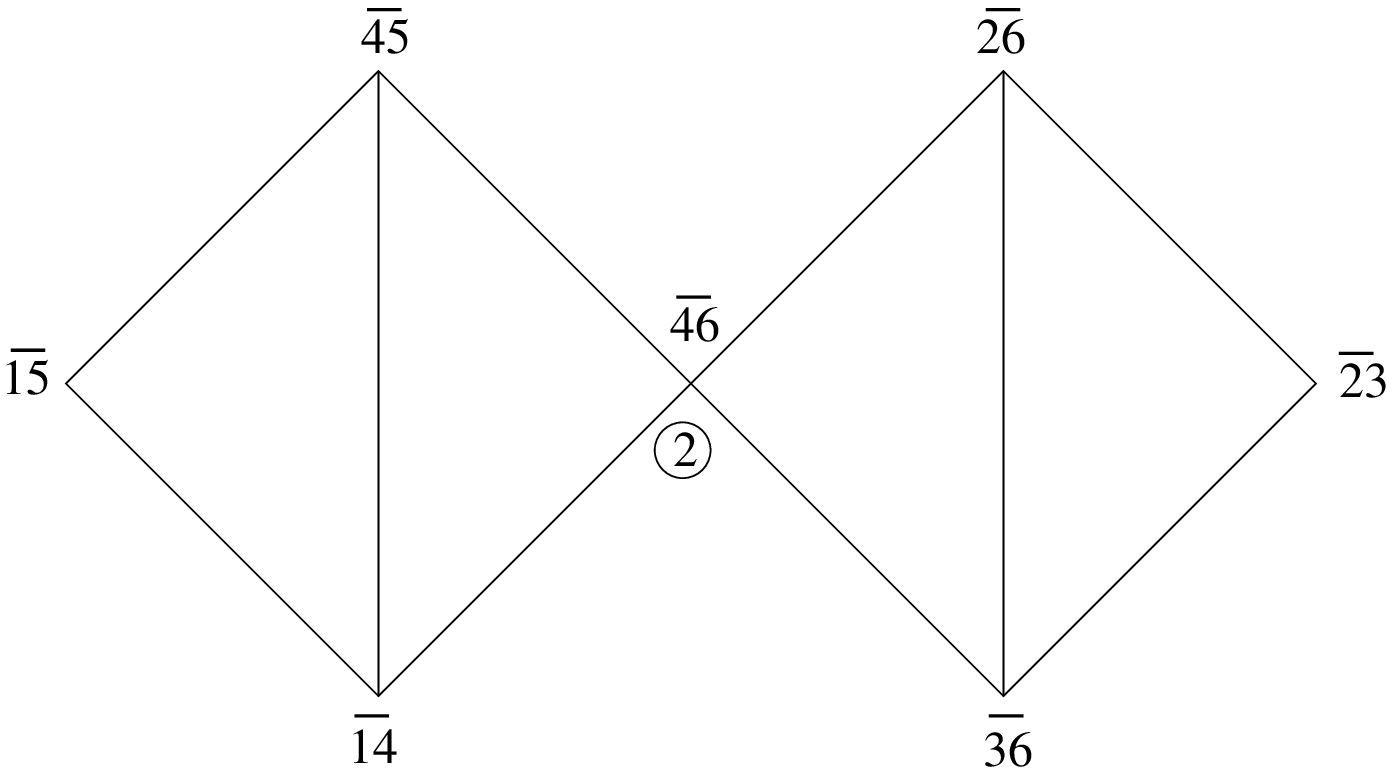}
\end{center}
\begin{center}
Figure 4 : $ R(u_4), R(u_5) $ in the cone $NHM_6^3 $
\end{center}

The cone $P_6^3$ has 4065 facets divided into at least 11 orbits. Below we
give representatives $f_1,\dots,f_{11}$ of the orbits $F_1,\dots,F_{11}$.
Their $(adjacency,incidence)$ pairs are, respectively, (1526,49), (703,41),
(100,23), (37,19), (31,18), (30,18), (23,17), (23,15), (22,18),
(18,16), (14,14). The facets $f_1, f_2$ are nonnegativity and
4-simplex facets; $f_{11}$ is a simplex cone. The $R$-graphs of the facets
$f_3$ and $f_4$ are on Figure 5.

\begin{center}
\epsfysize=4cm
\epsffile{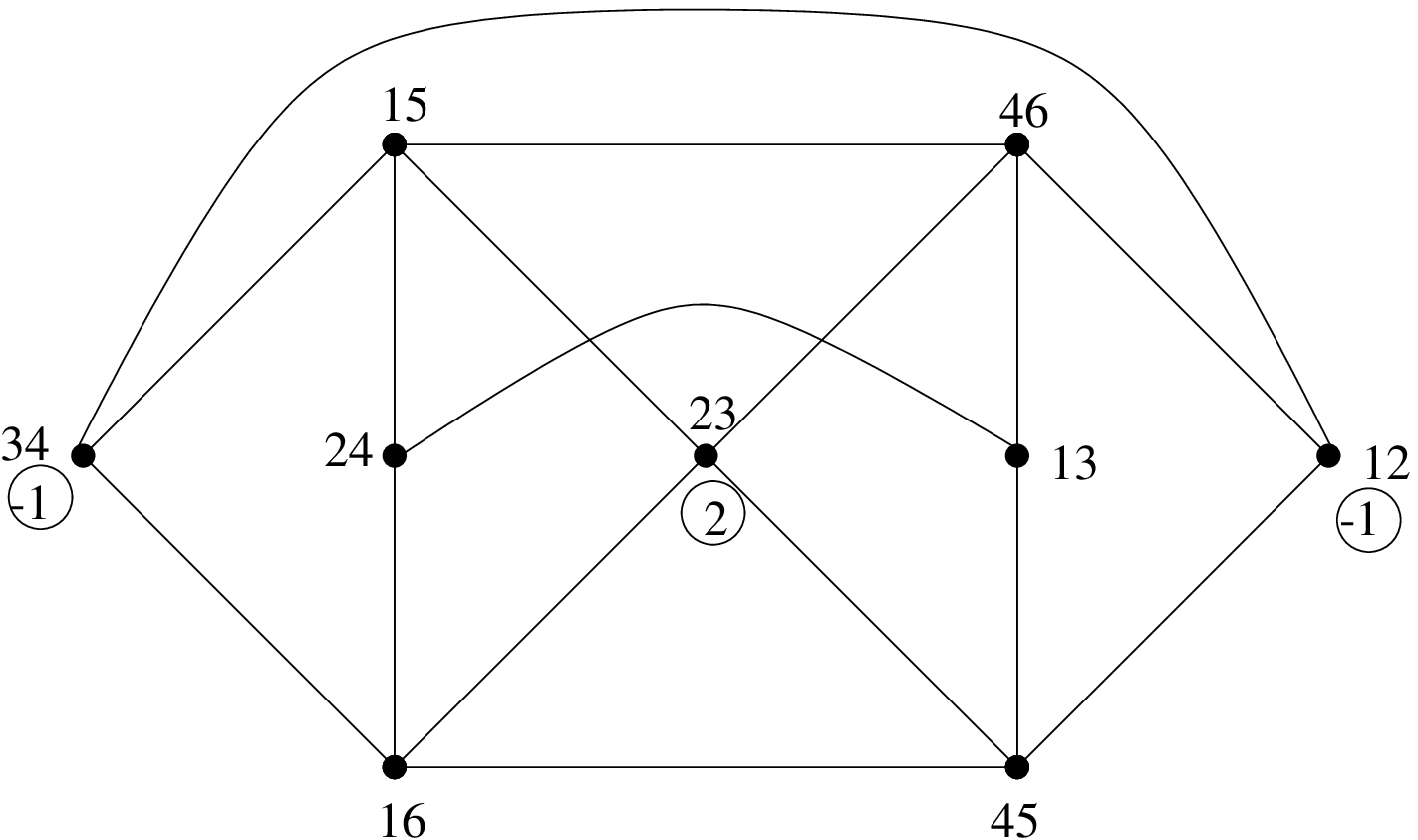}
\hspace{1cm}
\epsfysize=4cm
\epsffile{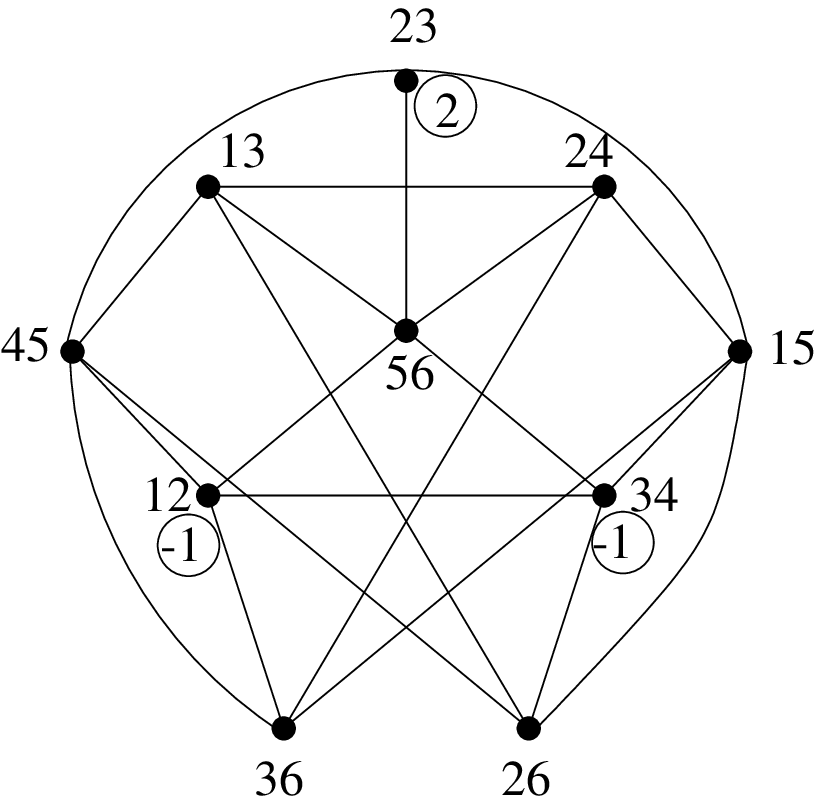}
\end{center}
\begin{center}
Figure 5 : $ \overline{R(f_3)}, \overline{R(f_4)} $ in the cone $P_6^3 $
\end{center}

$f_1= (0, 0, 0, 0, 0, 1, 0, 0, 0, 0, 0, 0, 0, 0, 0)$;

$f_2= (0, 1, 0, 0, 0, 1, 0, 0, 0,-1, 1, 1, 0, 0, 0)$;

$f_3=(-1, 1, 0, 1, 1, 2, 1, 0, 0,-1, 0, 0, 1, 1, 0)$;

$f_4=(-1, 1, 0, 1, 0, 2, 1, 0, 1,-1, 0, 1, 1, 0, 1)$;

$f_5=(-1, 1, 1, 2, 1, 2, 2,-1, 0,-2, 1, 0, 1, 2, 1)$;

$f_6=(-1, 1, 0, 2, 2, 2, 1,-1, 1,-1, 1,-1, 2, 2, 0)$;

$f_7=(-1, 1, 1, 3, 2, 2, 2,-2, 1,-2, 2,-1, 2, 3, 1)$;

$f_8= (1,-1, 3, 1, 4, 2, 2,-2, 1,-2, 2, 3, 2,-1, 1)$;

$f_9=(-1, 1, 1, 2, 2, 2, 2,-1,-1,-2, 1, 1, 1, 1, 2)$;

$f_{10}=(-1, 1, 1, 1, 2, 1, 1, 1,-1,-1,-1, 1, 2, 1, 1)$;

$f_{11}=(-1, 1, 1, 2, 0, 2, 2,-1, 1,-2, 1, 1, 1, 1, 2)$.

 \begin{table}
\caption{The adjacencies of extreme rays in the cone $NHM_6^3$}
 \begin{tabular}{|c|c|c|c|c|c|c|c|c|c|} \hline Orbit & Representative & 
 $O_{1}$ & $O_{2}$ & $O_{3}$ & $O_{4}$ & $O_{5}$ &Adj.&Inc. & $\vert 
 O_{i} \vert $ \\ \hline
 $O_{1}$ &$u_1=\alpha(1,2,3,456)$&$19$&$36$&$36$&$18$& $27$ & $136$ & $21, 12$ & $20$ \\
 $O_{2}$ &$u_2=\alpha(1,2,34,56)$&$16$&$18$&  $24$&$20$&$16$&$94$&$18, 11$&$45$\\
 $O_{3}$ &$u_3$ &$10$&$15$& $20$ & $15$ & $10$ & $70$ & $15, 10$ & $72$ \\
 $O_{4}$ &$u_4$ & $6$&$15$ & $18$&$9$ & $6$ & $54$ & $12, 9$ & $60$ \\
 $O_{5}$ &$u_5$ & $6$&$8$&$8$ & $4$ & $0$ & $26$ & $10, 8$& $90$ \\
  \hline 
\end{tabular}
\label{tab:tabl18}
 \end{table}

\begin{prop}

 The diameter of the skeleton graph of $P_6^3$ is $2$. Moreover:

(i) $G(O_1) = K_{20}$,  $G(O_2) = K_{45} - 15K_3$;

(ii) all non-edges are represented by $ \alpha(12,34,5,6)$ that are
nonadjacent to

$ \alpha(12,3,4,56), \alpha(1,2,34,56)$ (from the same
orbit $O_2$) and to

$ \alpha(125,3,4,6),\alpha(126,3,4,5), \alpha(345,1,2,6), \alpha(346,1,2,5)  $.

\end{prop}

In fact, both non-neighbors of $ \alpha(12,34,5,6) $ are in $O_2$. For both types
of non-edges - $ \alpha(12,34,5,6) $ with $ \alpha(12,3,4,56)$ and 
$ \alpha(125,3,4,6) $ - the ray $ \alpha(13,24,5,6)$ is a common neighbor.
Also, all 9 non-neighbors of a ray from  $O_1$, form $K_9$ in the skeleton
graph.

Notice that the skeleton of  $P_6^3$ is not an induced subgraph of the
skeleton of $NHM_6^3$; the only difference is in their restriction  $G(O_2)$
to the orbit of rays, represented by $u_2$. One can check that all neighbors
of a partition hemimetric  $ \alpha(a_1,a_2,b_1b_2,c_1c_2)$ from the same
orbit $O_2$ of $NHM_6^3$ are the 10 rays obtained by a transposition $(xy)$ and
the 8 rays obtained by a product $(a_1b_i)(a_2c_j)$ or $(a_1c_i)(a_2b_j)$ of two
transpositions. But in the skeleton of $P_6^3$, the ray
 $ \alpha(a_1,a_2,b_1b_2,c_1c_2)$ is adjacent to all other members of  $O_2$,
 except for the two rays, obtained from it by $(a_1b_1)(a_2b_2)$ or
$(a_1c_1)(a_2c_2)$. The complement of the graph, induced by all 18 neighbors
of the ray $ \alpha(a_1,a_2,b_1b_2,c_1c_2)$ from the same orbit $O_2$ of
$NHM_6^3$, is $C_4+C_4$ on 8 rays, obtained by a product $(a_1b_i)(a_2c_j)$
or a product $(a_1c_i)(a_2b_j)$, the skeleton of the cube on 8 rays, obtained by
$(a_ib_j)$ or $(a_ic_j)$, and it is $ \overline {K_2}$ on two rays obtained
by $(b_ic_j)$.

 \begin{table}
\caption{The adjacencies of extreme rays in the cone $P_6^3$}
 \begin{tabular}{|c|c|c|c|c|c|c|} 
 \hline 
 Orbit & Representative & $O_{1}$ &  $O_{2}$ & Adj.&Inc. & $\vert O_{i} 
 \vert $ \\ 
 \hline 
 $O_{1}$ & $\alpha(1,2,3,456)$ & $19$ & $36$ & $55$&$1113$ & $20$ \\ 
 $O_{2}$ & $\alpha(1,2,34,56)$ & $16$ & $42$ & $58$&$993$  & $45$ \\ 
 \hline 
 \end{tabular}
\label{tab:tabl9}
\end{table}

The skeletons of  $P_6^3$ and $NHM_6^3$ both contain a dominating clique
$O_1$; so their diameters are 2 or 3. In order to see closer the skeleton of
$NHM_6^3$, we now describe the local graph, denoted by $H$, of the ray
$u_5$. All 26 neighbors are in orbits $O_1,O_2,O_3,O_4$ only. It will be easier
to describe $ \overline H $. The restrictions
of $ \overline H$ on them are $ \overline {K_6}$, $C_8$, the skeleton of the
cube and $2K_2$, respectively. Two vertices from $O_1$ (say $15$ and $16$) are
isolated; so the diameter of $H$ is 2. Here we denote by $ij$ the $j$-th
member of the orbit $O_i$ in $H$. All edges of $ \overline H$
(without isolated vertices $15$ and $16$) are presented on
Figure 6. On the right picture the members of  $O_1$ are excluded while on the
left one the  members of  $O_2$ are excluded. $ \overline H$ does not
contain cross-edges among orbits $O_1$ and $O_2$.

\begin{center}
\epsfysize=4cm
\epsffile{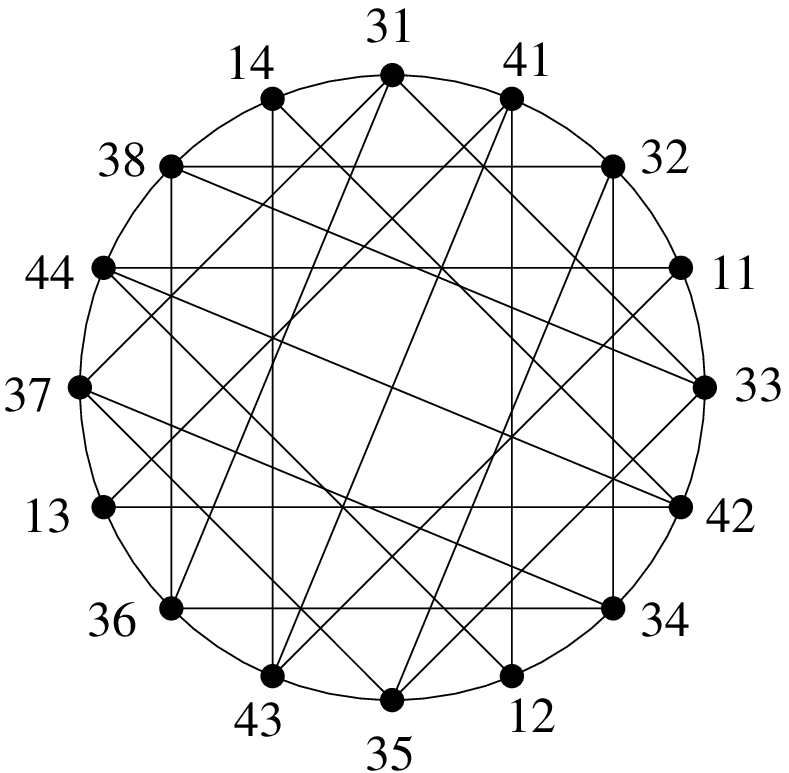}
\hspace{1cm}
\epsfysize=4cm
\epsffile{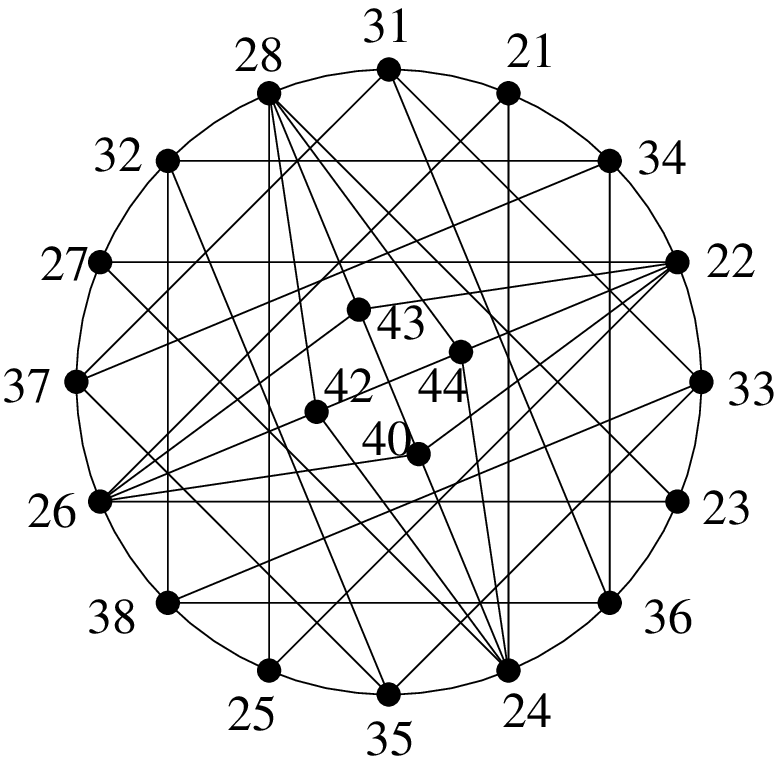}
\end{center}
\begin{center}
Figure 6 : A presentation of the local graph of a ray of the orbit $
O_5 $ of the cone $ NHM_6^3 $
\end{center}

\begin{prop} The ridge graph of $NHM_6^3$ has diameter $2$. Moreover:

(i) any 4-simplex facet $T_{ijkl,m}$ is adjacent to all but 5 facets:
$N_{ijkl}$ and all 4 other 4-simplex facets with the same support;  

(ii) the restrictions of the ridge graph to the orbits $F_1$ and $F_2$ are

$K_{5,5,5,5,5,5}$ and $K_{15}$, respectively.

\end{prop}

\begin{table}
\caption{The adjacencies of facets in the cone $NHM_6^3$}
 \begin{tabular}{|c|c|c|c|c|c|c|} 
 \hline 
 Orbit & Representative & $F_{1}$ &$F_{2}$& Adj.&Inc. & $\vert F_{i} \vert$ \\
\hline
 $F_{1}$ & $T_{1234,5}$ & $25$ & $14$ & $39$&$131$  & $30$ \\ 
 $F_{2}$ & $N_{1234}$ & $28$ & $14$ & $42$&$181$  & $15$ \\ 
 \hline 
 \end{tabular}
\label{tab:tabl10}
\end{table}

\section{Small 4-hemimetrics}

The cone $NHM_7^4$ has 3692 extreme rays divided into 8 orbits. We give below 
representatives $w_1,\dots,w_8$ of their orbits $O_1,\dots,O_8$. 
These vectors are indexed by 5-subsets of the set $\{1,\dots,7\}$; the 5-subsets are given as the complements of 2-subsets. 
The $(adjacency, incidence)$ pairs of those rays are, respectively, (985,48), 
(535,43), (315,38), (192,33), (126,28), (67,30), (43,25), (42,25). 
The first five vectors are $(0,1)$-valued; their graphs $R$  are $C_3, C_4, C_5, C_6, C_7 $, respectively.
The first two are partition $4$-hemimetrics; they represent both orbits
of $P_7^4$. The vectors $w_i, 1 \le i \le 4,$ and $w_6$ have same $R$-graphs
as the
members of orbits $O_i, 1 \le i \le 5,$ of $NHM_6^3$, respectively; so, the
graphs of Figure 5 represent also $w_4$ and $w_6$. The graphs
$R(w_7)$ and $R(w_8)$ are on Figure 7.

  $(\overline{12}, \overline{13}, \overline{14}, \overline{15}, \overline{16},
 \overline{17}, \overline{23}, \overline{24}, \overline{25}, \overline{26},
 \overline{27}, \overline{34}, \overline{35}, \overline{36}, \overline{37},
\overline{45}, \overline{46}, \overline{47}, \overline{56}, 
\overline{57}, \overline{67})$:

$w_1=(0, 0, 0, 1, 1, 0, 0, 0, 0, 0, 0, 0, 0, 0, 0, 0, 0, 0, 1, 0, 0)$;

$w_2=(0, 0, 0, 1, 1, 0, 0, 0, 1, 1, 0, 0, 0, 0, 0, 0, 0, 0, 0, 0, 0)$;

$w_3=(0, 0, 0, 1, 1, 0, 1, 0, 0, 1, 0, 0, 1, 0, 0, 0, 0, 0, 0, 0, 0)$;

$w_4=(0, 0, 0, 1, 1, 0, 1, 0, 0, 1, 0, 1, 0, 0, 0, 1, 0, 0, 0, 0, 0)$;

$w_5=(0, 0, 0, 1, 1, 0, 1, 0, 0, 1, 0, 1, 0, 0, 0, 0, 0, 1, 0, 1, 0)$;

$w_6=(0, 0, 0, 1, 1, 0, 0, 0, 0, 0, 0, 1, 0, 2, 1, 0, 0, 1, 1, 0, 0)$;

$w_7=(0, 0, 0, 2, 2, 0, 0, 0, 0, 1, 1, 1, 1, 0, 0, 1, 0, 0, 0, 0, 1)$;

$w_8=(0, 0, 0, 1, 1, 0, 0, 1, 0, 0, 1, 2, 1, 1, 0, 0, 0, 1, 0, 0, 0)$.

\begin{center}
\epsfysize=2.4cm
\epsffile{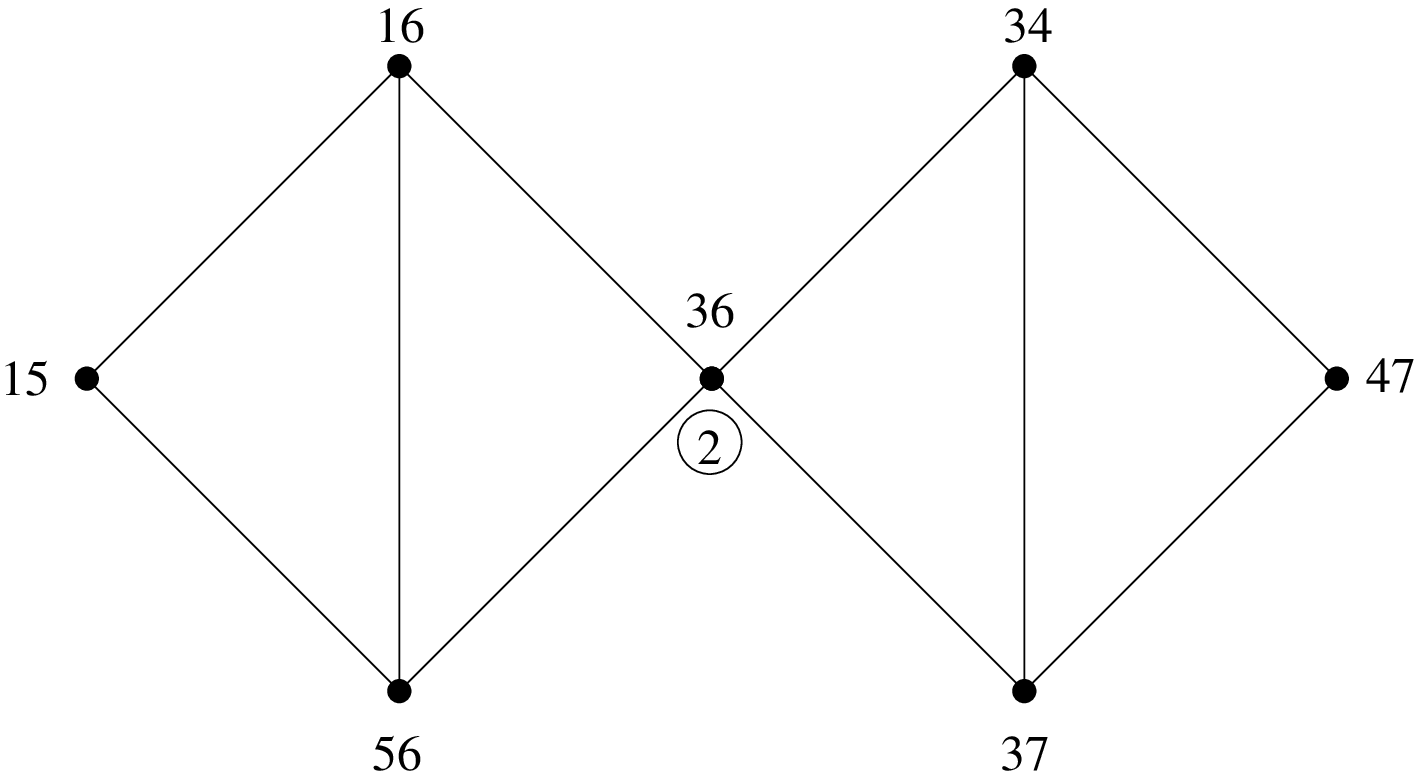}
\hspace{0.3cm}
\epsfysize=2.4cm
\epsffile{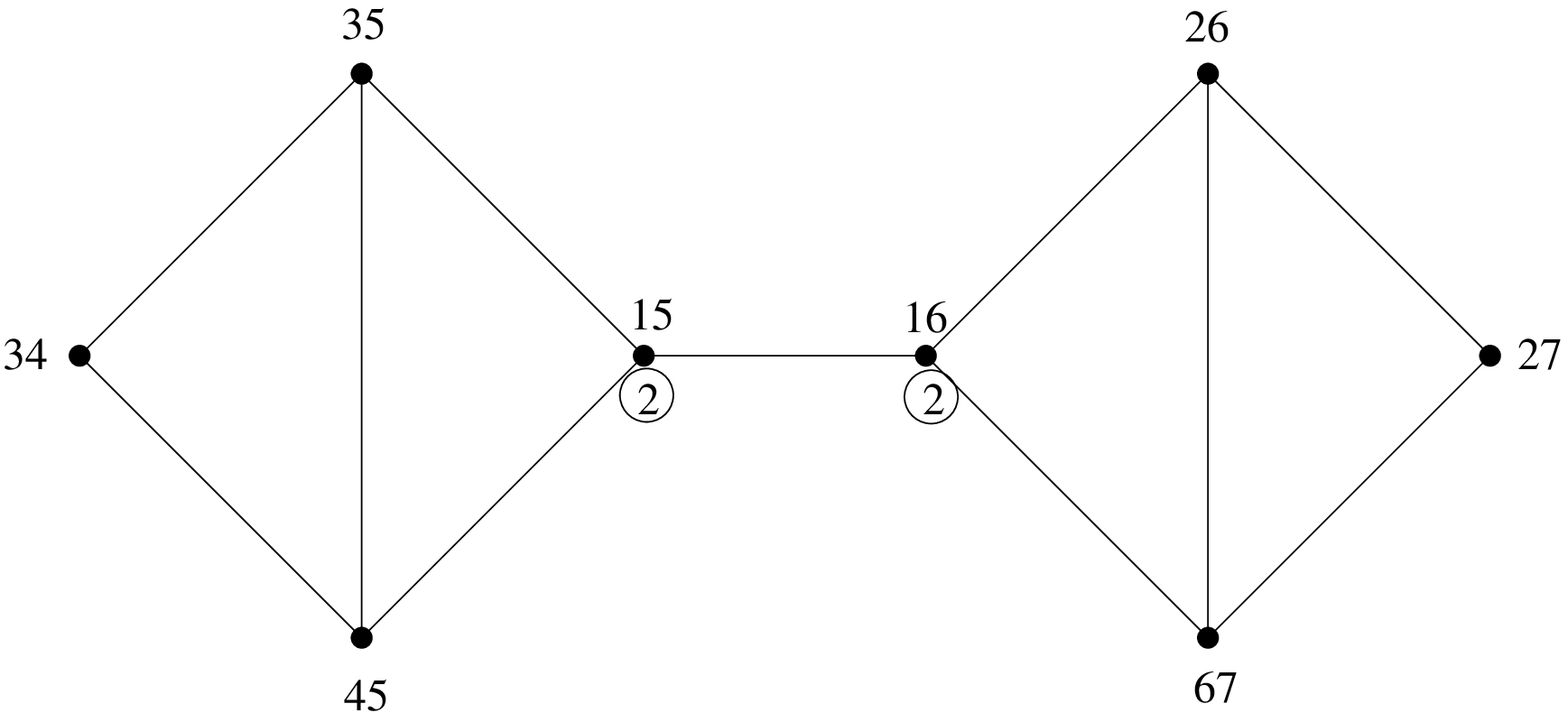}
\hspace{0.3cm}
\epsfysize=2.4cm
\epsffile{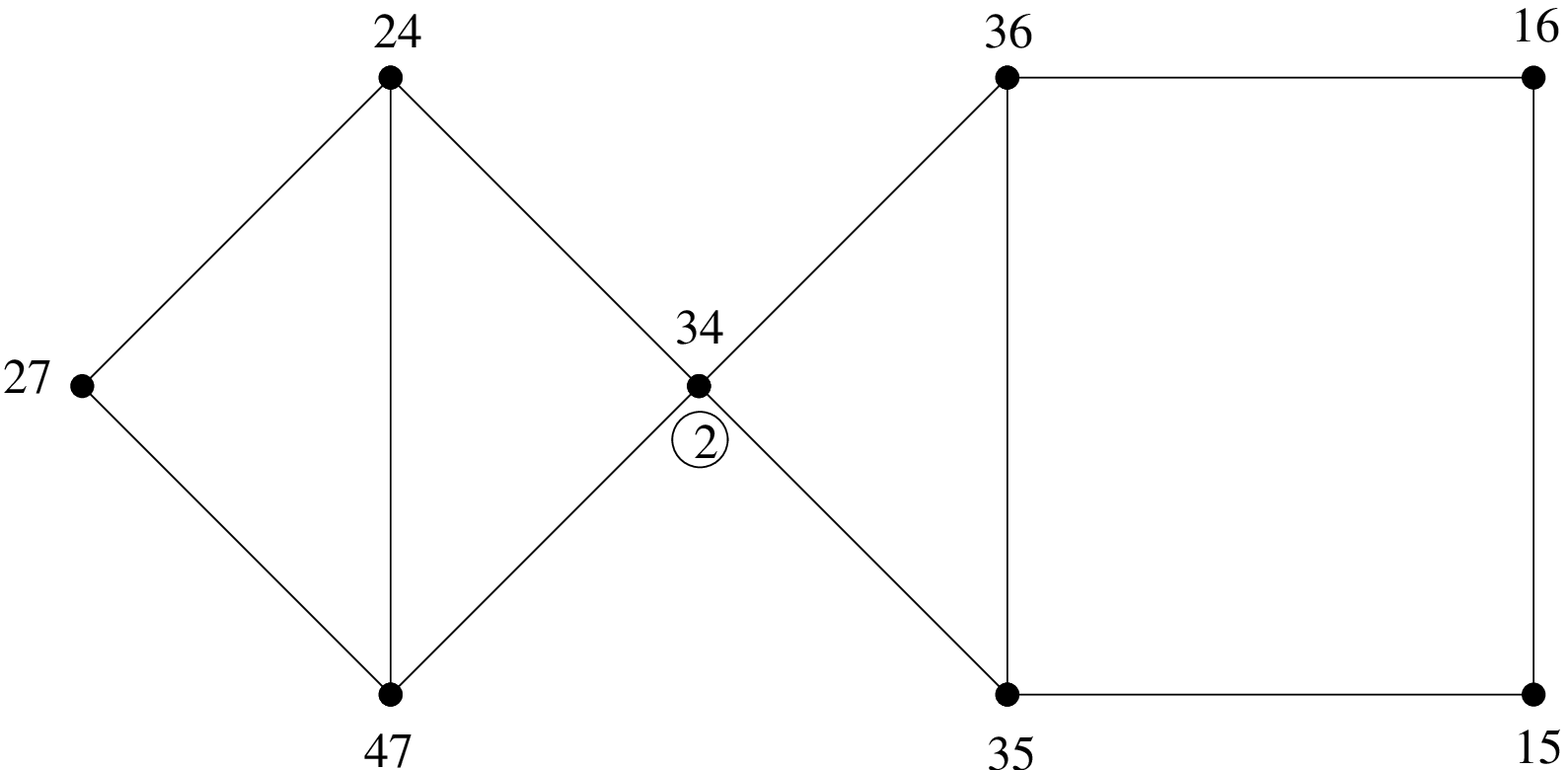}
\end{center}
\begin{center}
Figure 7 : $ R(W_6), R(W_7), R(W_8) $ in the cone $NHM_7^4 $
\end{center}

\begin{table}
\caption{The adjacencies of facets in the cone $NHM_7^4$}
 \begin{tabular}{|c|c|c|c|c|c|c|} 
 \hline 
 Orbit & Representative & $F_{1}$ &$F_{2}$& Adj.&Inc. & $\vert F_{i} \vert$ \\
\hline
 $F_{1}$ & $T_{12345,6}$ & $36$ & $20$ & $56$&$1302$  & $42$ \\ 
 $F_{2}$ & $N_{12345}$ & $40$ & $20$ & $60$&$2437$  & $21$ \\ 
 \hline 
 \end{tabular}
\label{tab:tabl11}
\end{table}

It is easy to check that the ridge graph of $NHM_7^4$ is $K_{6,6,6,6,6,6,6}$
on $F_1$ and $K_{21}$ on $F_2$. All non-edges among $F_1, F_2$ are of the form
$T_{i_1 \dots i_5,i_6}$ and  $N_{i_1 \dots i_5}$.

\section{Comparison of the small cones}

Now we compare some semimetric and m-hemimetric cones on $n$ points 
 for small  $n$. The triangle inequalities suffice to 
 describe the cut cones for $n \leq 4$,  but $CUT_{n} \subset  MET_{n}$ 
(strictly) for $n \geq 5$. 
The complete description of all the facets of the cut cone $CUT_{n}$ 
is known for $n \leq 8$, the complete description 
 of the semimetric cone $MET_{n}$ is known for 
 $n \leq 7$ (see, for example, the linear description of $MET_{7}$ in
~\cite{G}). Here the ``combinatorial explosion'' starts from $n = 8$.
 The number of orbits of 
facets and of extreme rays of those and other cones, when it
is known, is given in Table \ref{tab:tabl12}.

In fact, $P_n^2=NHM_n^2$ holds only for the smallest value $n=4$.
For $n=4, 5$ we  computed all facets, extreme rays and their adjacencies and 
incidences for three cones $P_n^2, NHM_n^2, HM_n^2$. For 2-hemimetrics the ``combinatorial
 explosion''(in terms of the amount of computation and  memory) starts already
for the cone $P_6^2$.
 
In the Table \ref{tab:tabl12} we compare the small 2-hemimetric cones
$P_n^2, NHM_n^2$ with the 1-hemimetric cones $CUT_n, MET_n$ and their
generalization in another direction: the cones $OMCUT_n, QMET_n$. Last
two cones consist of all {\it quasi-semimetrics} on $V_n$ and of
those obtained from {\it oriented multicuts}; see \cite{DP} for the notions
and results for them given in the Table \ref{tab:tabl12}.
The cones $NHM_n^2$ and $QMET_n$ have, besides of generalizations of the
usual
triangle inequality, only nonnegativity facets. In the Table \ref{tab:tabl12},
 columns 3 and 4 give the number of extreme rays and facets, respectively; in
parenthesis are given the numbers of their orbits. In column 5
are given the diameters of the skeleton and the ridge graphs of the cone 
specified in the row.
In the Table \ref{tab:tabl12}, the number of orbits of extreme rays and the diameter 
for cones
$QMET_5$, $NHM_6^2$, $P^2_6$, $P^4_7$ and dual $P_6^3$, $OMCUT_5$
are taken from recent work \cite{DuDe}, as well as the exact value of the diameter for $NHM_6^3$, $NHM_7^4$ and for the duals of $P_5^2$, $Cut_7$.

Incidences (to the extreme rays) of facets $T_{ijk,l}$ and  $N_{ijk}$ on the
cones $P_4^2=NHM_4^2$, $P_5^2, NHM_5^2$ amount to  3, 14 and 22, respectively,
but they are different (4001 and 3939) on $HM_6^2$. Incidences of similar
facets $T_{ij,k}$ ({\it oriented triangular inequality}, i.e.
$d(x,y)\leq d(x,z)+d(z,y)$ for a quasimetric $d$
), $N_{ij}$ (nonnegativity inequality)
are equal (to 7, 43) on cones $OMCUT_3=QMET_3, OMCUT_4$, but they are different
(78 and 80) on $QMET_4$.

For $n=4, 5$ we observe that the ridge graphs of $HM_n^2$ and $NHM_n^2$ are
induced subgraphs of the ridge graphs of $NHM_n^2$ and $P_n^2$, respectively.
The similar property does {\it not} hold for the 1-skeletons of those cones.
For example, any extreme ray of the orbit $O_2$ is adjacent to 14,6,2 members
of the same orbit in the cones $P_n^2, NHM_n^2, HM_n^2$, respectively.
Also, the ridge graph of $QMET_4$ is an induced subgraph of the ridge graph
of $OMCUT_4$, but the skeleton of $OMCUT_4$ is not an induced subgraph of the
skeleton of $QMET_4$ (see \cite{DP}). On the other hand, the ridge graph of
$MET_n$ and the skeleton of $CUT_n$ (for any $n$) have
diameters 2 and 1, respectively, and those graphs are induced subgraphs of the
ridge graph of $CUT_n$ and of the skeleton of $MET_n$, respectively (see Lemma
2.1 and Theorem 3.5 in \cite{DD1})

 \begin{table}
\caption{Some parameters of cones for small $n$}
 \scriptsize
\begin{tabular}{|c|c|c|c|c|} 
 \hline 
 cone &dimension &ext. rays (orbits) &facets (orbits)& diameters \\ 
 \hline 
 $P_{m+2}^m$=$NHM_{m+2}^m$& m+2& ${m+2} \choose 2$
(1)&$2m+4$ (2) & 2; 2 \\
     $m \ge 3$&&& &   \\
 \hline
  $CUT_3=MET_3$ &3&3(1)&3(1) & 1;  1 \\
 $P_4^2$=$NHM_4^2$ &4&6(1)&8(2) & 2;  3 \\
 $CUT_4=MET_4$ &6&7(2)&12(1) & 1;  2 \\ 
 $OMCUT_3=QMET_3$ &6&12(2)&12(2) & 2;  2 \\
 $CUT_5$ &10&15(2)&40(2) & 1;  2 \\ 
 $MET_5$ &10&25(3)&30(1) & 2;  2 \\  
 $P_5^2$ &10&25(2)&120(4) & 2;  3 \\ 
 $NHM_5^2$ &10&37(3)&30(2) & 2;  2\\ 
 $OMCUT_4$ &12&74(5)&72(4) & 2;  2 \\ 
 $QMET_4$ &12&164(10)&36(2) & 3;  2 \\
 $CUT_6$ &15&31(3)&210(4) & 1;  3 \\
 $MET_6$ &15&296(7)&60(1) & 2;  2 \\
 $P_6^3$ &15&65(2)&4065(16) & 2;  3 \\
 $NHM_6^3$ &15&287(5)&45(2) & 3;  2 \\
 $P_6^2$ &20&90(3)&$ \ge 2095154( \ge 3086)$ & 2;  ? \\ 
 $NHM_6^2$ &20&12492(41)&80(2) & 3;  2 \\ 
 $OMCUT_5$ &20&540(10)&35320(194) & 2;  3 \\ 
 $QMET_5$ &20&43590(229)&80(2) & 3;  2 \\
 $P_7^4$ &21&140(2)& 474390(153) & 2; 3 \\  
 $NHM_7^4$ &21&3692(8)&63(2) & 3;  2 \\ 
 $CUT_7$ &21&63(3)&38780(36) & 1;  3 \\ 
 $MET_7$ &21&55226(46)&105(1) & 3;  2 \\ 
 $CUT_8$ &28&127(4)&$ \ge 49604520( \ge 2169)$ & 1;  ? \\
 $P_8^5$ &28&266(2)&$ \ge 322416108( \ge 8792)$ & ?;  ? \\
 $NHM_8^5$&28&55898(13)&84(2) & 3;  2 \\
\hline 
 \end{tabular}
\label{tab:tabl12}
\end{table}

 \section{ Conjectures for general $m, n$}
  
\begin{conj}
The two partition $m$-hemimetrics
$\alpha(S_1,\dots,S_{m+1})$ and $\alpha(T_1,\dots,T_{m+1})$ on $V_n$ are nonadjacent in the
skeleton of $P_n^m$ if and only if there exist six different subsets
$S_i,S_j,S_k$ and $T_{i'},T_{j'},T_{k'}$, such that $S_i \bigcup S_j =T_{k'}$
and $S_k=T_{i'} \bigcup T_{j'}$.

The conjecture holds for $m=1$: all cut semimetrics are adjacent. It holds for
$n-m=2$:
we have the graph $J(m+2,2)$. It also holds for $(m,n)=(2,5)$ and $(3,6)$. 
\end{conj}
\begin{conj}
 The ridge graphs of $HM_n^m$ 
  and of $NHM_n^m$ are induced subgraphs of the ridge graphs of
 $NHM_n^m$ and $P_n^m$, respectively.

Recall that the ridge graph of $NHM_n^m$ has two orbits of
vertices: $F_1, F_2$, consisting of $(n-m-1)$$n \choose {m+1}$ simplex and
$n \choose {m+1}$ nonnegativity inequalities.
\end{conj}
\begin{conj}

The ridge graph $NHM_n^m$ satisfies:
  
   (i) The $(m+1)$-simplex facet $T_{i_1 \dots i_{m+1},i_{m+2}}$ is adjacent
to all other facets, except the following $m+2$ facets:

all other $(m+1)$-simplex facets with the same support and $N_{i_1 \dots i_{m+1}}$; 

   (ii) $G(F_2)= \overline {J(n,3)}$ for $m=2$ and
$G(F_2)=K_{n \choose {m+1}}$ for $m \ge 3$.

Clearly, (i) implies that the restriction of the ridge graph on $F_1$ is
 $G(F_1)=K_{{m+2},\dots,{m+2}}$. It is easy to see that Conjecture $3$ would imply that the diameter of the
ridge graph of $NHM_n^m$ is 2 (it was proved in \cite{DD1} that the diameter
of the ridge graph of $NHM_n^1 = MET_n$ is 2). In fact, to see it for
$m=2$ consider all 3 types of
pairs of nonadjacent vertices:

(i) let $x,y \in F_1$ have the same support, say, ${1234}$. Suppose that
$x_{124}=y_{124}=-1$. Then $N_{123}$ is a common neighbor for $x$ and $y$.

(ii) for $N_{123}$ and $N_{124}$, any tetrahedron facet
$T_{134,2}$ is their common neighbor.

(iii) for $N_{123}$ and $T_{123,4}$, the facet $N_{345}$ is a common neighbor.
\end{conj}
\begin{conj}

The extreme rays of $NHM_n^m$ include:

(i) any   ray whose $R$-graph is an $R$-graph of an extreme ray of
 $NHM_{n-1}^{m-1}$;

(ii) every $(0,1)$-valued extreme ray of $NHM_{m+3}^m$ with
$R$-graph $C_i$ ($3 \le i \le {m+3}$).
\end{conj}

\end{document}